\documentclass[12pt]{amsart}
\usepackage{amssymb,amsmath,amsthm}
\usepackage{calrsfs}
\usepackage{graphicx}
\newtheorem{lm}{Lemma}[section]

\newtheorem{pr}[lm]{Proposition}
\newtheorem{thm}[lm]{Theorem}
\newtheorem{ls}[lm]{Corollary}

\newtheorem{prop}{Property}
\begin{document}
\large
\title{Non-commutative $L_p$-spaces associated with a Maharam trace}
\author{Vladimir Chilin}
\author{Botir Zakirov}
\address{Vladimir Chilin \\Department of Mathematics, National University of Uzbekistan\\
Vuzgorodok, 100174 Tashkent, Uzbekistan} \email{{\tt
chilin@ucd.uz}}
\address{Botir Zakirov \\Tashkent Railway Engineering Institute \\
Odilhodjaev str. 1, 100167 Tashkent, Uzbekistan} \email{{\tt
botirzakirov@list.ru}}
\date{}

\begin{abstract}
Non-commutative $L_p$-spaces $L^p(M,\Phi)$ associated with the Maharam trace are defined and their  dual spaces are described.\\



\end{abstract}

\maketitle

\section{Introduction}

Development of the theory of integration for measures $\mu$ with the values in Dedekind complete Riesz spaces
has inspired the study of $(bo)$-complete lattice-normed spaces $L^p(\mu)$ (see, for example, \cite{Ku1}, 6.1.8).
Note that, if the measure $\mu$ satisfies the Maharam property, then the spaces $L^p(\mu)$ are Banach-Kantorovich.

The existence of center-valued traces on finite von Neumann algebras naturally leads to a study of
the integration for traces with the values in a complex Dedekind complete Riesz space $F_\mathbb{C}=F\oplus iF.$
For commutative von Neumann algebras, the development of $F_\mathbb{C}$-valued integration is a
part of the study of the properties of order continuous positive maps of Riesz spaces,
for which we refer to the treatise by A.G. Kusraev  \cite{Ku1}.
The operators possessing the Maharam property provide important examples of such mappings, while the
$L^p$-spaces associated with such operators are non-trivial examples of Banach-Kantorovich Riesz spaces.

Let $M$ be a non-commutative von Neumann algebra,
$F_\mathbb{C}$ a von Neumann subalgebra in the center of $M$,
and let $\Phi: M\to F_\mathbb{C}$ be a trace such that $\Phi(zx)=z\Phi(x)$ for all $z\in F_\mathbb{C},~x\in M.$
Then the non-commutative $L^p$-space $L^p(M,\Phi)$ is
a Banach-Kantorovich space \cite{Chil.G.1.}, \cite{Chil.K.}, and the trace
$\Phi$ satisfies the Maharam property, that is, if $0\leq z\leq
\Phi(x),~z\in F_\mathbb{C},~0\leq x\in M,$ then there exists
$0\leq y \leq x$ such that $\Phi(y)=z $ (compare with \cite{Ku1},
3.4.1).

In \cite{ChZ}, a faithful normal trace $\Phi$ on $M$ with the values in an arbitrary complex Dedekind complete Riesz space was considered.
In particular, a complete description of such traces in the case when $\Phi$ is a Maharam trace was given.
In the same paper, utilizing the locally measure topology on the algebra $S(M)$ of all measurable operators affiliated with $M,$
the Banach-Kantorovich space $L^1(M,\Phi) \subset S(M)$ was constructed and a version of Radon-Nikodym-type theorem for Maharam traces was established.

In the present article, we define a new class of Banach-Kantorovich spaces, non-commutative
$L_p$-spaces $L^p(M,\Phi)$ associated with a Maharam trace; also, we give a description of their dual spaces.
We use the terminology and results of the theory of von Neumann algebras (\cite{SZ.}, \cite{Take1.}),
the theory of measurable operators (\cite{Seg.}, \cite{Mur_m}), and of the theory of Dedekind
complete Riesz space and Banach-Kantorovich spaces (\cite{Ku1}).

\section{Preliminaries}

Let $X$ be a vector space over the field $\mathbb C$ of complex numbers, and let $F$ be a Riesz space. A mapping $\|\cdot \|: X\to F$ is said to be a vector ($F$-valued) norm if it satisfies the following axioms:
\begin{enumerate}
\item  $\| x\|\ge 0,$ $\| x\|=0$ $\Leftrightarrow$ $x=0$ ($x\in X$);
\item $\| \lambda x\|=|\lambda|\ \| x\|$ ($\lambda\in \mathbb C,$ $x\in X$);
\item $\| x+y\| \le \| x\| +\| y\|$ ($x,y\in X$).
\end{enumerate}

A norm $\|\cdot\|$ is called decomposable   if the following property holds:

\begin{prop}\label{decomposability} If $f_1,f_2\ge 0$ and $\|x\|=f_1+f_2,$ then there exist $x_1, x_2\in X$ such that $x=x_1+x_2$ and $\|x_k\|= f_k$ ($k=1,2$).
\end{prop}

If property \ref{decomposability} is valid only for disjoint elements $f_1, f_2\in F,$ the norm is called  disjointly decomposable or, briefly,  d-decomposable.

The pair $(X,\|\cdot\|)$ is called  a lattice-normed space (shortly, LNS). If the norm $\| \cdot\|$ is decomposable (d-decomposable), then so is the space $(X,\|\cdot\|).$

A net $\{x_\alpha\}_{\alpha\in A}\subset X$ $(bo)$- converges to $x\in X$ if the net $\{\| x_\alpha -x\|\}_{\alpha\in A}$ $(o)$-converges to zero in the Riesz space $F.$ A net $\{x_\alpha \}_{\alpha \in A}$ is said to be a $(bo)$- Cauchy net if $\sup\limits_{\alpha,\beta\ge \gamma} \| x_\alpha-x_\beta\| \downarrow 0.$ An LNS is called $(bo)$- complete if any $(bo)$-Cauchy net $(bo)$-converges.  A Banach-Kantorovich space (shortly, BKS) is a  d-decomposable $(bo)$-complete LNS. It is well known that every BKS is a decomposable LNS.

Let $F$ be a Dedekind complete Riesz space, and let  $F_\mathbb{C}=F\oplus iF$ be the complexification of $F.$ If $z=\alpha + i\beta \in F_\mathbb{C},~\alpha, \beta \in F, $ then $\overline{z}:=\alpha -i\beta,$ and $|z|:=\sup\{Re( e^{i\theta}z): 0\leq \theta < 2\pi \}$ (see\cite{Ku1}, 1.3.13).

Let $(X,\|\cdot\|_X)$ be the BKS over $F.$ A linear operator $T:X\to F_\mathbb{C}$ is said to be $F$-bounded if there exists $0\leq c\in F$ such that $|T(x)|\leq c\|x\|_X$ for all $x\in X.$ For any $F$-bounded operator $T,$  define the element $\|T\|=\sup\{|T(x)|:x\in X,$ $\|x\|_X\leq \mathbf{1}_F\},$ which is called the abstract $F$-norm of the operator $T$ (\cite{Ku1}, 4.1.3). It is known that $|T(x)|\leq \|T\|\,\|x\|_X$ for all $x\in X$ (\cite{Ku1}, 4.1.1).

The set $X^*$ of all $F$-bounded linear mappings from $X$ into $F_\mathbb{C}$ is called the $F$-dual space to the BKS $X.$ For $T,S\in X^*,$ we set $(T+S)(x)=Tx+Sx,$ $(\lambda T)(x)=\lambda Tx,$ where $x\in X,~\lambda\in \mathbb{C}.$ It is clear that  $X^*$ is a  linear space  with respect to the introduced algebraic operations. Moreover, $(X^*, \|\cdot\|)$ is a BKS (\cite{Ku1}, 4.2.6).

Let $H$ be a Hilbert space, let $B(H)$ be the $*$-algebra of all bounded linear operators on $H,$ and let $\mathbf{1}$ be the identity operator on $H.$  Given  a von Neumann algebra $M$ acting on $H,$ denote by $Z(M)$  the center of $M$ and by $P(M)$  the lattice of all projections in $M$. Let  $P_{fin}(M)$ be the set of all finite projections in $M.$

A densely-defined  closed linear operator $x$ (possibly unbounded) affiliated with $M$ is said to be \emph{measurable} if there exists a sequence $\{p_n\}_{n=1}^{\infty}\subset P(M)$ such that $p_n\uparrow \mathbf{1}$, \ $p_n(H)\subset \mathfrak{D}(x)$ and $p_n^\bot=\mathbf{1}-p_n \in P_{fin}(M) $ for every $n=1,2,\ldots$ (here $\mathfrak{D}(x)$ is the domain of $x$). Let us denote by $S(M)$ the set  of all measurable  operators.

Let $x,y$ be measurable  operators. Then $x+y,~xy$ and $x^*$ are densely-defined and preclosed. Moreover, the closures $\overline{x+y}$ (strong sum), $\overline{xy}$ (strong product) and $x^*$ are also measurable, and $S(M)$ is a  $*$-algebra with respect to the strong sum, strong product, and the adjoint operation (see \cite{Seg.}). For any subset $E\subset S(M)$ we denote by $E_h$ (resp. $ E_+$ ) the set of all self-adjoint (resp. positive ) operators from $E.$

For  $x\in S(M)$  let  $x=u|x|$ be the polar decomposition, where $|x|=(x^*x)^{\frac{1}{2}},$ $u$ is a partial isometry in $B(H).$ Then $u\in M$ and $|x|\in S(M).$ If $x\in S_h(M)$ and $\{E_\lambda(x)\}$ are the spectral  projections of $x,$ then $\{E_\lambda(x)\}\subset P(M).$

Let $M$ be a commutative von Neumann algebra. Then  $M$ is $*$-isomorphic to the $*$-algebra $L^\infty(\Omega,\Sigma,\mu)$ of all essentially   bounded complex measurable functions  with the identification almost everywhere, where $(\Omega,\Sigma,\mu)$ is a measurable space.   In addition $S(M)\cong L^0(\Omega,\Sigma,\mu),$ where $L^0(\Omega,\Sigma,\mu)$ is the $*$-algebra  of all complex measurable functions  with the identification almost everywhere \cite{Seg.}.

The locally measure topology  $t(M)$ on $L^0(\Omega,\Sigma,\mu)$ is by definition the linear (Hausdorff) topology whose fundamental system of neighborhoods of $0$ is given by $$ W(B,\varepsilon,\delta)=\{f\in\ L^0(\Omega,\, \Sigma,\, \mu) \colon \hbox{ there exists a set } \ E\in \Sigma, \mbox{ such that} $$ $$  \ E\subseteq B, \ \mu(B\setminus E)\leqslant\delta, \ f\chi_E \in L^\infty(\Omega,\Sigma,\mu), \ \|f\chi_E\|_{{L_\infty}(\Omega,\Sigma,\mu)}\leqslant\varepsilon\}. $$ Here \ $\varepsilon, \ \delta $ run over all strictly positive numbers and  $B\in\Sigma$, \ $\mu(B)<\infty.$ It is known that $(S(M),t(M))$ is a complete topological $*$-algebra.

It is clear that zero neighborhoods $W(B,\varepsilon,\delta)$ are closed and have the following  property: if $f\in W(B,\varepsilon, \delta),\, g\in L^\infty(\Omega,\Sigma,\mu), \|g\|_{L{\infty}(\Omega,\Sigma,\mu)}\leq 1,$ then $gf\in W(B,\varepsilon,\delta).$

A net $\{f_\alpha\}$ converges  locally in measure to $f$ (notation: $f_\alpha \stackrel{t(M)}{\longrightarrow}f$) if and only if  $f_\alpha \chi_B $ converges in $\mu$-measure to $f\chi_B$   for each $B\in\Sigma$ with $\mu(B)<\infty$.   If  $M$ is $\sigma$-finite  then there exists a faithful finite normal  trace $\tau$ on $M.$ In this case, the topology $t(M)$ is metrizable, and convergence  $f_n \stackrel{t(M)}{ \longrightarrow} f$ is equivalent to convergence in  trace $\tau$ of the sequence $f_n$ to $f.$

Let now $M$ be an arbitrary finite von Neumann algebra, $\Phi_M: M\to Z(M)$ be a center-valued trace on $M$ (\cite{SZ.}, 7.11). Let $Z(M)\cong L^\infty(\Omega,\Sigma,\mu).$ The locally measure topology  $t(M)$ on $S(M)$ is  the linear (Hausdorff) topology whose fundamental system of neighborhoods of $0$ is given by $$ V(B,\varepsilon, \delta ) = \{x\in S(M)\colon \ \mbox{there exists } \ p\in P(M), z\in P(Z(M)) $$  $$ \mbox{ such that} \ xp\in M, \|xp\|_{M}\leqslant\varepsilon, \ z^\bot \in W(B,\varepsilon,\delta),  \ \Phi_M(zp^\bot)\leqslant\varepsilon z\},$$ where $\|\cdot\|_{M}$ is the $C^*$-norm in $M.$  It is known that  $(S(M),t(M))$ is a complete topological $*$-algebra \cite{Yead1}.

From (\cite{Mur_m}, \S 3.5) we have the following criterion for convergence  in the topology $t(M). $

\begin{pr}\label{2.1.} A net $\{x_\alpha\}_{\alpha\in A} \subset S(M)$ converges to zero in the topology $t(M)$   if and only if $\Phi_M(E^\bot_\lambda (|x_\alpha|) \stackrel{t(M)}{\longrightarrow} 0$ for any $\lambda>0.$
\end{pr}

Let $M$ be an arbitrary  von Neumann algebra, and let $F$ be a Dedekind complete Riesz space. An {\em $F_\mathbb{C}$-valued trace} on the von Neumann algebra $M$ is a linear mapping $\Phi:M\to F_\mathbb{C}$ with $\Phi(x^*x)=\Phi(xx^*)\geq 0$ for all $x\in M.$ It is clear that $\Phi(M_h)\subset F, ~\Phi(M_+)\subset F_+=\{a\in F: a\geq 0\}.$ A trace $\Phi$ is said to be {\em faithful} if the equality $\Phi(x^*x)=0$ implies $x=0,$ {\em normal} if $\Phi(x_\alpha)\uparrow\Phi(x)$ for every $x_\alpha,x\in M_h,~x_\alpha \uparrow x.$

If $M$ is a finite von Neumann algebra, then  its canonical center-valued trace $\Phi_M:M\to Z(M)$ is an example of a $Z(M)$-valued faithful normal trace.

Let us list some properties of the trace $\Phi:M\to F_\mathbb{C}.$

\begin{pr}(\cite{ChZ}) $(i)$ Let $x,y,a,b \in M.$ Then

$\Phi(x^*)=\overline{\Phi(x)},$ $ \Phi(xy)=\Phi(yx),$ $\Phi(|x^*|)=\Phi(|x|),$

$|\Phi(axb)| \leq \|a\|_M \|b\|_M \Phi(|x|);$

$(ii)$ If $ \Phi$ is a faithful trace, then $M$ is  finite;

$(iii)$ If $x_n, x \in M$ and $\|x_n-x\|_M \to 0,$ then $|\Phi(x_n)-\Phi(x)|$ relative uniform converges  to zero;

$(iv)$ $\Phi(|x+y|)\leq \Phi(|x|)+\Phi(|y|)$ for all $x, y \in M.$
\end{pr}

The trace $\Phi:M\to F_\mathbb{C}$ possesses  {\em the Maharam property} if for any $x\in M_+,~0\leq f\leq \Phi(x), ~f\in F,$ there exists $y\in M_+,$  $y\leq x$ such that $\Phi(y)=f.$ A faithful normal $F_\mathbb{C}$-valued trace $\Phi$ with the Maharam property is called {\em a Maharam trace} (compare with \cite{Ku1}, III, 3.4.1). Obviously, any faithful finite numerical trace on $M$ is a $\mathbb{C}$-valued Maharam trace.

Let us give another examples of Maharam traces. Let $M$ be a finite von Neumann algebra, let $\mathcal{A}$ be a von Neumann subalgebra in $Z(M),$ and let $T:Z(M) \to \mathcal{A}$ be an injective linear positive normal operator. If $f\in S(\mathcal{A})$ is a reversible positive element, then $\Phi(T,f)(x)= f T(\Phi_M(x))$ is an $S(\mathcal{A})$-valued faithful normal trace on $M.$ In addition, if $T(ab)=aT(b)$ for all $a\in \mathcal{A}, b\in Z(M),$ then $\Phi(T,f)$ is a Maharam trace on $M.$

If $\tau$ is a faithful normal finite numerical trace on $M$ and $ \dim(Z(M))>1,$ then $\Phi(x)=\tau(x) \mathbf{1}$ is a $Z(M)$-valued faithful normal trace, which does not possess the Maharam property (see \cite{ChZ}).

Let $F$ have a weak order unit $\mathbf{1}_F.$ Denote by $B(F)$ the complete Boolean algebra of unitary elements with respect to $\mathbf{1}_F,$ and let $Q$ be the Stone compact space of the Boolean algebra $B(F).$ Let $C_\infty(Q)$ be the Dedekind complete Riesz space of all continuous functions $a: Q\to [-\infty, +\infty]$ such that $a^{-1}(\{\pm \infty\})$ is a  nowhere dense subset of $Q.$ We identify $F$ with the order-dense ideal in $C_\infty(Q)$ containing  algebra $C(Q)$ of all continuous real functions on $Q.$ In addition, $\mathbf{1}_F$ is identified with the function equal to 1  identically on $Q$ (\cite{Ku1}, 1.4.4).

We need the following theorem from \cite{ChZ}.

\begin{thm} \label{2.2.} Let $\Phi$ be an $F_\mathbb{C}$-valued Maharam trace on a von Neumann algebra $M.$
Then there exists a von Neumann subalgebra $\mathcal{A}$ in $Z(M),$ a $*$-isomorphism $\psi$ from $\mathcal{A}$ onto
the $*$-algebra $C(Q)_\mathbb{C},$ a positive linear  normal operator $\mathcal{E}$ from
$Z(M)$ onto $\mathcal{A}$ with $\mathcal{E}(\mathbf{1})=\mathbf{1},~\mathcal{E}^2=\mathcal{E},$ such that

$1)$ $\Phi(x)=\Phi(\mathbf{1})\psi(\mathcal{E}(\Phi_M(x)))$ for all $x \in M;$

$2)$ $\Phi(zy)=\Phi(z\mathcal{E}(y))$ for all $z,y \in Z(M);$

$3)$ $\Phi(zy)=\psi(z)\Phi(y)$ for all $z\in \mathcal{A},$ $y\in M.$
\end{thm}

Due to Theorem \ref{2.2.},  the $*$-algebra $\mathcal{B}=C(Q)_\mathbb{C}$  is a commutative von Neumann algebra, and  $*$-algebra $C_\infty(Q)_\mathbb{C}$ is identified with the  $*$-algebra $S(\mathcal{B}).$ It is clear that the $*$-isomorphism $\psi$ from $\mathcal{A}$ onto $\mathcal{B}$ can be extended to a $*$-isomorphism from $S(\mathcal{A})$ onto $S(\mathcal{B}).$ We denote this mapping also by $\psi.$

Let $\Phi$ be a $S(\mathcal{B})$-valued Maharam trace  on a von Neumann algebra $M.$   A net $\{x_\alpha\}\subset S(M)$ converges to $x\in S(M)$ with respect to the trace $\Phi$ (notation: $x_\alpha \stackrel{\Phi}{\longrightarrow}x$) if $\Phi(E_\lambda^\bot(| x_\alpha-x|))\stackrel{t(\mathcal{B})}{\longrightarrow}0$ for all $\lambda>0.$

\begin{pr}\label{2.3.}(\cite{ChZ}) $x_\alpha \stackrel{\Phi}{\longrightarrow}x$ iff $x_\alpha \stackrel{t(M)}{\longrightarrow}x.$
\end{pr}

An operator  $x\in S(M)$ is said to be {\em $\Phi$-integrable} if there exists a sequence $\{x_n\} \subset M$ such that $x_n \stackrel{\Phi}{\to} x $ and $\|x_n-x_m\|_\Phi \stackrel{t(\mathcal{B})}{\longrightarrow} 0$ as $n,m \to \infty.$

Let $x$ be a $\Phi$-integrable operator from $ S(M).$ Then  there exists a $\widehat{\Phi}(x)\in S(\mathcal{B})$ such that $\Phi(x_n) \stackrel{t(\mathcal{B})}{\longrightarrow}\widehat{\Phi}(x).$ In addition $\widehat{\Phi}(x)$ does not depend on the choice of a sequence $\{x_n\}\subset M,$ for which $x_n\stackrel{\Phi}{\longrightarrow}x,$ $\Phi(|x_n-x_m|) \stackrel{t(\mathcal{B})}{\longrightarrow}0$ \cite{ChZ}. It is clear that each operator $x\in M$ is $\Phi$-integrable and $\widehat{\Phi} (x)=\Phi(x).$

Denote by $L^1(M,\Phi)$ the set of all $\Phi$-integrable operators from $S(M).$ If $x\in S(M)$ then $x\in L^1(M,\Phi)$ iff $|x|\in L^1(M,\Phi),$  in addition $|\widehat{\Phi}(x)|\leq \widehat{\Phi}(|x|)$ \cite{ChZ}. For any $x\in L^1(M,\Phi),$ set $\|x\|_{1,\Phi}=\widehat{\Phi}(|x|).$ It is known that $L^1(M,\Phi)$ is a linear subspace of $S(M),$ $ML^1(M,\Phi)M \subset L^1(M,\Phi),$ and $x^*\in L^1(M,\Phi)$ for all $x\in L^1(M,\Phi)$ \cite{ChZ}. Moreover,  the following theorem   is true.

\begin{thm}\label{2.4.}(\cite{ChZ}) $(i)$ $(L^1(M,\Phi), \| \cdot\|_{1,\Phi})$ is a Banach-Kantorovich space;

$(ii)$ $S(\mathcal{A})L^1(M,\Phi) \subset L^1(M,\Phi),$ in addition $\widehat{\Phi}(zx)=\psi(z)\widehat{\Phi}(x)$ for all $z\in S(\mathcal{A}),~x\in L^1(M,\Phi).$
\end{thm}

\section{$L_p$-spaces associated with a Maharam  trace}

Let $\mathcal{B}$ be a commutative von Neumann algebra, which is
$*$-isomorphic to a von Neumann subalgebra $\mathcal{A}$ in  $Z(M),$  and let $\Phi:M\to
S(\mathcal{B})$ be a Maharam trace on $M$ (see Theorem
\ref{2.2.}). For any $p>1,$ set $L^p(M,\Phi)=\{x\in S(M): |x|^p\in
L^1(M,\Phi)\}$ and $\|x\|_{p,\Phi}=\widehat{\Phi}(|x|^p)
^{\frac{1}{p}}.$ It is clear that $M\subset L^p(M,\Phi).$

Let $e$ be a nonzero projection in $\mathcal{B},$ and put
$\Phi_e(a)=\Phi(a)e,~a\in M.$ A mapping $\Phi_e: M\to
S(\mathcal{B}e)$ is a normal (not necessarily faithful) $S(\mathcal{B}e)$-valued trace on $M.$
Denote by $s(\Phi_e):=\mathbf{1}-\sup\{p\in P(M): \Phi_e(p)=0\}$ the support
of the trace $\Phi_e.$ It is clear that $s(\Phi_e)\in P(Z(M))$ and $\Phi_e(a)=\Phi(as(\Phi_e))$  is a faithful normal $S(\mathcal{B}e)$-valued trace on $Ms(\Phi_e)$ (compare \cite{SZ.}, 5.15). Moreover $\Phi_e$ possesses   the Maharam property.

If $e$ and $g$ are  orthogonal nonzero projections in
$P(\mathcal{B}),$ then $\Phi_g(s(\Phi_e))=\Phi(s(\Phi_e))g=
\Phi_e(\mathbf{1})g=\Phi(\mathbf{1})eg=0,$ i.e.
$s(\Phi_e)s(\Phi_g)=0.$ Let $\{e_i\}_{i\in I}$ be a family of nonzero mutually orthogonal projections in $P(\mathcal{B})$ with
$\sup\limits_{i\in I}e_i=\mathbf{1}_\mathcal{B},$ where $\mathbf{1}
_\mathcal{B}$ is the unit of the algebra $\mathcal{B}.$ If $z=\mathbf{1}-\sup\limits_{i\in I}s(\Phi_{e_i})$ then $\Phi(z)e_i=\Phi_{e_i}(z)=0$ for all $i\in I.$ Therefore
$\Phi(z)=0,$ i.e. $z=0,$ or $\sup\limits_{i\in I}s(\Phi_{e_i})=\mathbf{1}.$

Further, we need the following

\begin{pr}\label{3.1.} Let $ x\in S(M)$ and let $\{e_i\}_{i\in I}$
be the family of nonzero mutually orthogonal projections in
$P(\mathcal{B})$ with $\sup_{i\in I}e_i=\mathbf{1}_\mathcal{B}.$
Then $x\in L^p(M,\Phi)$ if and only if $xs(\Phi_e)\in
L^p(Ms(\Phi_{e_i}),\Phi_{e_i})$ for all $i\in I.$ In addition
$\|x\|_{p,\Phi}e_i=\|xs(\Phi_{e_i})\|_{p,\Phi_{e_i}}.$
\end{pr}
\begin{proof} Let $x\in L^p(M,\Phi),$ $a_n=E_n(|x|^p)|x|^p$
where $E_n(|x|^p)$ is the spectral projection of $|x|^p$
corresponding to the interval $(-\infty, n].$ It is clear that
$a_n\stackrel{\Phi}{\longrightarrow}|x|^p$ and
$\Phi(|a_n-a_m|)\stackrel{t(\mathcal{B})}{\longrightarrow}0$ as
$n,m \to \infty.$ Hence, $a_ns(\Phi_{e_i}) \stackrel{\Phi_{e_i
}}{\longrightarrow} |x|^ps(\Phi_{e_i})$ (see Proposition
\ref{2.3.}). In addition, from the inequality $\Phi_{e_i}
(|a_ns(\Phi_{e_i})-a_ms(\Phi_{e_i})|)=\Phi(|a_n-a_m|s(\Phi_{e_i}))\leq
\Phi(|a_n-a_m|),$ we have $\Phi_{e_i}(|a_ns(\Phi_{e_i})-
a_ms(\Phi_{e_i})| \stackrel{t(\mathcal{B}e_i)}{\longrightarrow}0.$
This means that $|xs(\Phi_{e_i})|^p=|x|^ps(\Phi_{e_i})\in
L^1(Ms(\Phi_{e_i}),\Phi_{e_i})$ and $\|xs(\Phi_{e_i})\|_{p,
\Phi_{e_i}}= \widehat{\Phi}_{e_i}(|x|^ps(\Phi_{e_i}))^{
\frac{1}{p}}=(\widehat{\Phi}(|x|^p)e_i)^{\frac{1}{p}}=\|x\|_{p,\Phi}e_i.$

Conversely, let $xs(\Phi_{e_i})\in L^p(Ms(\Phi_{e_i}),\Phi_{e_i})$
for all $i\in I.$ Set $a_{n,i}=E_n(|xs(\Phi_{e_i})|^p) |xs(
\Phi_{e_i})|^p.$ It is clear that $a_{n,i}\uparrow |xs(\Phi_{e_i})|^p=
|x|^ps(\Phi_{e_i})$ as $n\to \infty$ for any fixed $i\in I.$
Therefore $a_{n,i}\stackrel{t(Ms(\Phi_{e_i}))}{\longrightarrow}
|x|^ps(\Phi_{e_i}),$ $\Phi_{e_i}(|a_{n,i}-a_{m,i}|) \stackrel{t(
\mathcal{B}e_i)}{\longrightarrow}0$ as $n,m \to \infty.$  Since
$0\leq \Phi(\sqrt{a_{n,i}}a_{m,j}\sqrt{a_{n,i}})=
\Phi(a_{n,i}a_{m,j})\leq \|a_{m,j}\|_M\Phi(a_{n,i})=\|a_{m,j}\|_
M\Phi (a_{n,i})e_i$ and $\Phi(a_{n,i}a_{m,j})\leq
\|a_{n,i}\|_M\Phi(a_{m,j})e_j,$ we have $\Phi(a_{n,i}a_{m,j})=0.$ Hence,
$a_{n,i}a_{m,j}=0$ for all $n,m,~i\neq j.$ Since $0\leq
a_{n,i}\leq ns(\Phi_{e_i}),$ $s(\Phi_{e_i})s(\Phi_{e_j})=0,$
$i\neq j,$  there is an $x_n\in M_+$ such that
$x_ns(\Phi_{e_i})=a_{n,i}.$ Using the equality $\sup\limits_{i\in
I}s(\Phi_{e_i})=\mathbf{1},$ we obtain $x_n\stackrel{t(M)
}{\longrightarrow} |x|^p$ (\cite{Zak}), moreover $\Phi(|x_n-x_m|)
\stackrel{t(\mathcal{B})}{\longrightarrow}0.$ Therefore $x\in
L^p(M,\Phi).$
\end{proof}

Similar to in the case of the space $L^1(M,\Phi),$ the subset
$L^p(M,\Phi)$ is invariant with respect to the action of
involution in  $S(M).$ The following proposition is devoted to this fact.

\begin{pr}\label{3.2.} If $x\in L^p(M,\Phi),$ then $x^*\in L^p(M,\Phi)$ and $\|x\|_{p,\Phi}=\|x^*\|_{p,\Phi}.$
\end{pr}
\begin{proof}

Let $x=u|x|$ be the polar decomposition of $x.$ Since an algebra $M$ has a finite type, we can suppose that $u$ is a
unitary operator in $M.$  For each $y\in S(M),$ we set $U(y)=uyu^*.$ Then  the mapping $U:S(M)\to S(M)$
is a $*$-isomorphism, and therefore $U(\varphi(y))= \varphi(U(y))$ for any continuous function
$\varphi:[0,+\infty)\to [0,+\infty)$ and $y\in S_+(M)$ \cite{Zak}.
If $\varphi(t)=t^p,$ $p>1,$ $t\geq 0,$ and $y\in S_+(M)$ then
$uy^pu^*=(uyu^*)^p.$ In particular,  we obtain the equality $|x^*|^p=u|x|^pu^*.$ Hence,  $x^*\in L^p(M,\Phi).$ Moreover $\|x^*\|_{p,\Phi}=\widehat{\Phi}(|
x^*|^p)^{\frac{1}{p}} =\widehat{\Phi}(u|x|^pu ^*)^{\frac{1}{p}}
=\widehat{\Phi}(|x|^p)^{\frac{1}{p}}=\|x\|_{p,\Phi}.$
\end{proof}

Now we need a version of the H\"{o}lder inequality for Maharam traces. In the proof of this inequality for numerical traces,  properties of
decreasing rearrangements of integrable operators are used \cite{Yead2}.
For Maharam traces such theory of decreasing rearrangements does not exact.
Therefore  we use another approach connected with the
concept of a bitrace on a $C^*$-algebra.

Let $\mathcal{N}$ be a $C^*$-algebra. A function $s:\mathcal{N}
\times \mathcal{N} \to \mathbb{C}$ is called a bitrace on
$\mathcal{N}$ (\cite{Dix1}, 6.2.1) if  the following relations hold:

$(i)$ $s(x,y)$ is positively defined sesquilinear Hermitian form
on $\mathcal{N};$

$(ii)$ $s(x,y)=s(x^*,y^*)$ for all $x,y\in \mathcal{N};$

$(iii)$ $s(zx,y)=s(x,z^*y)$ for all $x,y,z\in \mathcal{N};$

$(iv)$ for any $z\in \mathcal{N},$ the mapping $x\to zx$  is continuous on $(\mathcal{N},
\|\cdot\|_s)$ where $\|x\|_s=\sqrt{s(x,x)},$ $x\in \mathcal{N};$

$(v)$ the set $\{xy: x,y\in \mathcal{N}\}$ is dense in
$(\mathcal{N},\|\cdot\|_s).$

If $\mathcal{N}$ has a unit, then condition $(v)$ holds
automatically.

Let us list examples of bitraces associated with the Maharam
trace.

Let $M$ be a von Neumann algebra, let $\Phi:M\to S(\mathcal{B})$ be a
Maharam trace and let $Q=Q(P(\mathcal{B}))$ be the Stone compact space of the
Boolean algebra $P(\mathcal{B}).$ We claim that
$s(\Phi(\mathbf{1}))=\mathbf{1}_\mathcal{B}.$ If it is not the case, then
$e=\mathbf{1}_\mathcal{B}-s(\Phi(\mathbf{1}))\neq 0$  and
$z=\psi^{-1}(e)\neq 0$ where $\psi$ is a $*$-isomorphism from
Theorem \ref{2.2.}. By Theorem \ref{2.4.}$(ii),$ we have
$\Phi(z)=e\Phi(\mathbf{1})=0,$ which contradicts to the faithfulness of
the trace $\Phi.$ Thus, $s(\Phi(\mathbf{1}))=\mathbf{1}_\mathcal{
B},$ and therefore the following elements are defined:
$(\Phi(\mathbf{1}))^{-1} \in S_+(\mathcal{B})$ and $(\Phi(
\mathbf{1}))^{-1}\Phi(x)\in C(Q)$ where $x\in M.$ For any $t\in
Q,$ set $\varphi_t(x)=(\Phi(\mathbf{1})^{-1}\Phi(x))(t).$ It is
clear that $\varphi_t$ is a finite numerical trace on $M.$ The function
$s_t(x,y)=\varphi_t(y^*x)=\varphi_t(xy^*)$ is a bitrace on $M.$ In
fact, the conditions $(i)-(iii)$ are obvious.  $(iv)$ follows from the inequality $\|zx\|_{s_t}=
\sqrt{\varphi_t((zx)^*(zx))} =\sqrt{\varphi_t(x^*z^*zx)}\leq
\|z\|_M\|x\|_{s_t}.$

Let $s(x,y)$ be  an arbitrary bitrace   on a von Neumann algebra
$M.$ Set $N_s=\{x\in M:s(x,x)=0\}.$ It follows from (\cite{Dix1},
6.2.2) that $N_s$ is a self-adjoint two-sided ideal in $M.$
We consider the factor-space $M/N_s$ with the scalar product $([x],
[y])_s=s(x,y)$ where $[x], [y]$ are the equivalence classes from
$M/N_s$ with representatives $x$ and $y,$ respectively. Denote by
$(H_s,(\cdot,\cdot)_s)$ the Hilbert space which is  the completion of
$(M/N_s, (\cdot,\cdot)_s).$ By the formula $\pi_s(x)([y])=[xy], ~x,y\in M,$ one defines a $*$-homomorphism $\pi_s:M\to B(H_s).$ In addition $\pi_s(\mathbf{1}_M) =\mathbf{1}_{B(H_s)}.$

Denote by $U_s(M)$ the von Neumann subalgebra in $B(H_s)$
generated by operators $\pi_s(x),$ i.e. $U_s(M)$ is the closure of
the $*$-subalgebra $\pi_s(M)$ in $B(H_s)$ with respect to the weak operator
topology. According to (\cite{Dix2}, s. 85-88), there exists  a faithful  normal semifinite numerical trace $\tau_s$ on
$(U_s(M))_+$ such that $\tau_s(\pi(x^2))=([x],[x])=s(x,x)$ for all $x\in M_+.$ If $\varphi$ is a trace  on $M$ and $s(x,y)=\varphi(y^*x)$ then $\tau_s(\pi_s(x^2))=
\varphi(x^2)$ for all $x\in M_+.$ This means that $\tau_s( \pi_s(x))
=\varphi(x)$ for any $ x\in M_+.$ In addition, if $\varphi
(\mathbf{1}_M)<\infty,$ then $\tau_s(\mathbf{1}_{B(H_s)})<\infty.$
Consequently, $\tau_s$ is a faithful normal finite trace on $U_s(M).$

\begin{thm}\label{3.3.} Let $\Phi$ be a $S(\mathcal{B})$-valued
Maharam trace on the von Neumann algebra $M,$ $p,q>1,$
$\frac{1}{p}+\frac{1}{q}=1.$ If $x\in L^p(M,\Phi),$ $y\in
L^q(M,\Phi),$ then $xy\in L^1(M,\Phi)$ and
$\|xy\|_{1,\Phi}\leq\|x\|_{p,\Phi}\|y\|_{q,\Phi}.$
\end{thm}
\begin{proof} We consider  the bitrace
$s_t(x,y)=\varphi_t(y^*x)$ on $M$ where $\varphi_t(x)=((\Phi( \mathbf{1})
)^{-1}\Phi(x))(t),~t\in Q(P(\mathcal{B})).$ Denote by $\tau_t$ a
faithful normal finite trace on $(U_{s_t}(M))_+$ such that $\tau_t(
\pi_{s_t}(x))=\varphi_t(x)$ for all $x\in M_+.$ Since the trace
$\tau_t$ is finite, $\tau_t(\pi_{s_t}(x)=\varphi_t(x)$ for any
$x\in M.$ Let $L^p(U_{s_t}(M), \tau_t) $ be the non-commutative $L^p$-space  associated with the numerical  trace
$\tau_t.$ It follows from \cite{Yead2} that
$$
\|\pi_{s_t}(x)\pi_{s_t}(y)\|_{1,\tau_t}\leq \| \pi_{s_t}
(x)\|_{p,\tau_t}\|\pi_{s_t}(y)\|_{q,\tau_t},
$$
i.e.
$$
\tau_t(|\pi_{s_t}(xy)|)\leq \tau_t(|\pi_{s_t}(x)|^p
)^{\frac{1}{p}}\tau_t(|\pi_{s_t}(y)|^q)^{\frac{1}{q}}.
$$
Since  $\pi_{s_t}(|x|)=|\pi_{s_t}(x)|,~x\in M,$ we get $\pi_{s_t}(|x|^p)=(\pi_{s_t}(|x|))^p$ (\cite{Dix1}, 1.5.3).

Thus,
$$
\tau_t(\pi_{s_t}(|xy|))\leq \tau_t(\pi_{s_t}(|x|^p))^{\frac{1}{p}}
\tau_t(\pi_{s_t}(|y|^q))^{\frac{1}{q}},$$ i.e.
$\varphi_t(|xy|)\leq \varphi_t(|x|^p)^{\frac{1}{p}}
\varphi_t(|y|^q)^{\frac{1}{q}},$ or
$$(\Phi(\mathbf{1}))^{-1}\Phi(|xy|)(t)\leq
[((\Phi(\mathbf{1}))^{-1}\Phi(|x|^p))(t)]^{\frac{1}{p}}[((\Phi(\mathbf{1}))^{-1}\Phi(|y|^q))(t)]^{\frac{1}{q}}$$
for all $t\in Q(P(\mathcal{B})).$ This means that
$$(\Phi(\mathbf{1}))^{-1}\Phi(|xy|)\leq
[((\Phi(\mathbf{1}))^{-1}\Phi(|x|^p))]^{\frac{1}{p}}[((\Phi(\mathbf{1}))^{-1}\Phi(|y|^q))]^{\frac{1}{q}}.$$
Multiplying this inequality by $\Phi(\mathbf{1}),$ we get
$\|xy\|_{1,\Phi}\leq\|x\|_{p,\Phi}\|y\|_{q,\Phi}.$

Let now $x\in L^p_+(M,\Phi),$ $y\in L^q_+(M,\Phi).$ We claim that
$xy\in L^1(M,\Phi).$ Set $a_n=E_n(x)x,~b_n=E_n(y)y.$ We have
$a_n,b_n\in M_+$ and $a_n\uparrow x,~b_n\uparrow y,$ in particular,
$a_n\stackrel{\Phi}{\longrightarrow}x,~
b_n\stackrel{\Phi}{\longrightarrow}y.$ Hence, $a_nb_n\in M$ and $a_nb_n\stackrel{\Phi}{\longrightarrow}xy.$ In addition,
$\|a_nb_n-a_mb_m\|_{1,\Phi}\leq\|a_nb_n-a_nb_m\|_{1,\Phi}+\|a_nb_m-a_mb_m\|_{1,\Phi}\leq\|a_n\|_{p,\Phi}\|b_n-b_m\|_{q,\Phi}+\|a_n-a_m\|_{p,\Phi}\|b_m\|_{q,\Phi}.$
Since $\|a_n\|_{p,\Phi}\leq \|x\|_{p,\Phi},~\|b_m\|_{q,\Phi}\leq
\|y\|_{q,\Phi},$ and for $n>m,$
$\|a_n-a_m\|^p_{p,\Phi}=\widehat{\Phi}(x^pE_n(x)E_m^\bot(x))\stackrel{t(\mathcal{B})}{\longrightarrow}0,$
$\|b_n-b_m\|^q_{q,\Phi}=\widehat{\Phi}(y^qE_n(y)E_m^\bot(y))\stackrel{t(\mathcal{B})}{\longrightarrow}0,$
we get
$\|a_nb_n-a_mb_m\|_{1,\Phi}\stackrel{t(\mathcal{B})}{\longrightarrow}0$
as $n,m\to \infty.$ This means that $xy\in L^1(M,\Phi)$ and
$\|a_nb_n-xy\|_{1,\Phi}\stackrel{t(\mathcal{B})}{\longrightarrow}0.$
The inequality
$|\|xy\|_{1,\Phi}-\|a_nb_n\|_{1,\Phi}|\leq\|xy-a_nb_n\|_{1,\Phi}$
implies
$\|a_nb_n\|_{1,\Phi}\stackrel{t(\mathcal{B})}{\longrightarrow}\|xy\|_{1,\Phi}.$
Since
$$\|a_nb_n\|_{1,\Phi}\leq\|a_n\|_{p,\Phi}\|b_n\|_{q,\Phi}\stackrel{t(\mathcal{B})}{\longrightarrow}\|x\|_{p,\Phi}\|y\|_{q,\Phi},$$ we obtain
$\|xy\|_{1,\Phi}\leq\|x\|_{p,\Phi}\|y\|_{q,\Phi}.$

If $x\in L^p(M,\Phi)$ is arbitrary, $y\in L^q_+(M,\Phi)$ and
$x=u|x|$ is the polar decomposition of $x$ with the unitary $u\in
M,$ then $xy=u(|x|y)\in L^1(M,\Phi)$ and $\|xy\|_{1,\Phi}=
\||x|y\|_{1,\Phi}\leq\|x\|_{p,\Phi}\|y\|_{q,\Phi}.$

Let now $x\in L^p(M,\Phi),$ $y\in L^q(M,\Phi)$ be arbitrary and let
$y^*=v|y^*|$ be the polar decomposition of $y^*$ with the unitary
$v\in M.$ According to Proposition \ref{3.2.}, $|y^*|\in
L^q(M,\Phi)$ and $\|y^*\|_{q,\Phi}=\|y\|_{q,\Phi}.$ Therefore
$xy=(x|y^*|)v^*\in L^1(M,\Phi)$ and
$$\|xy\|_{1,\Phi}=\|x|y^*|\|_{1,\Phi}
\leq\|x\|_{p,\Phi}\||y^*|\|_{q,\Phi}=\|x\|_{p,\Phi}\|y\|_{q,\Phi}.$$
\end{proof}

\begin{thm}\label{3.4.} Let $\Phi, M, p,$ and $q$ be the same as in Theorem
\ref{3.3.}. If $x\in S(M),$ $xy\in L^1(M,\Phi)$ for all $y\in
L^q(M,\Phi)$ and the set $D(x)=\{|\widehat{\Phi}(xy)|: y\in
L^q(M,\Phi),$ $\|y\|_{q,\Phi}\leq \mathbf{1}_\mathcal{B}\}$ is
bounded in $S_h(\mathcal{B}),$ then $x\in L^p(M,\Phi)$ and
$\|x\|_{p,\Phi} =\sup D(x).$
\end{thm}
\begin{proof} Let $x\neq 0,$ and let $x=u|x|$ be the polar decomposition of $x$ with the unitary
$u\in M.$ Set $y_n=|x|^{p-1}E_n(|x|)E^\bot_{\frac{1}{n}}
(|x|)u^*,~n=1,2,\dots$ It is clear that $y_n\in M$ and
$$xy_n=u|x|^pE_n(|x|)E^\bot_{\frac{1}{n}}(|x|)u^*=uE_n(|x|)E^\bot_{\frac{1}{n}}(|x|)|x|^pE_n(|x|)E^\bot_{\frac{1}{n}}(|x|)u^*\geq
0.$$ On the other hand,
$$|y_n|^2=uE_n(|x|)E^\bot_{\frac{1}{n}}(|x|)|x|^{2p-2}E_n(|x|)E^\bot_{\frac{1}{n}}(|x|)u^*=$$
$$uE_n(|x|)E^\bot_{\frac{1}{n}}(|x|)|x|^{\frac{2p}{q}}E_n(|x|)E^\bot_{\frac{1}{n}}(|x|)u^*,$$
and therefore $0\leq|y_n|^q=(|y_n|^2)^{\frac{q}{2}}=xy_n,$ in
particular, $\|y_n\|_{q,\Phi}=\Phi(xy_n)^{\frac{1}{q}}.$

Since $xy_n\stackrel{t(M)}{\longrightarrow}u|x|^pu^*\neq 0,$ we have
$xy_n\neq 0$ for all$n\geq n_0.$ Set $e_n=s(\Phi
(xy_n)$ as $n\geq n_0.$ Since $S_h(\mathcal{B})=C_\infty(Q(P(
\mathcal{B})),$  there exists a unique $b_n\in
S_+(\mathcal{B})e_n$ such that $b_n\Phi(xy_n)=e_n.$ It is clear that
$b_n^{\frac{1}{q}}\Phi^{\frac{1}{q}}(xy_n)=e_n.$ If
$z_n=\psi^{-1}(e_n),~a_n=\psi^{-1}(b_n^{\frac{1}{q}})\in
S(\mathcal{A}z_n),$ then by  theorem \ref{2.4.}$(ii),$
$a_ny_n\in L^q(M,\Phi)$ and $\|a_ny_n\|^q_{q,\Phi}=\widehat{\Phi}
(a_n^q|y_n|^q)=b_n\widehat{\Phi}(xy_n)=e_n\leq
\mathbf{1}_\mathcal{B}.$ Hence, $|\widehat{\Phi}(a_nxy_n)|=
|\widehat{\Phi}(x(a_ny_n))|\leq \sup D(x)$ for all $n\geq n_0.$ On
the other hand,
$$
\widehat{\Phi}(a_nxy_n)=b_n^{\frac{1}{q}}\widehat{\Phi}(xy_n)=
(b_n\widehat{\Phi}(xy_n))^{\frac{1}{q}}\widehat{\Phi}(xy_n)^{1-\frac{1}{q}}=
\widehat{\Phi} (xy_n)^{\frac{1}{p}}=
$$ $$
\widehat{\Phi}(u|x|^pE_n(|x|)E^\bot_{\frac{1}{n}}(|x|)u^*)^{\frac{1}{p}}=
\widehat{\Phi}(|x|^pE_n(|x|)E^\bot_{\frac{1}{n}}(|x|))^{\frac{1}{p}}.
$$
Since $(|x|^pE_n(|x|)E^\bot_{\frac{1}{n}}(|x|))\uparrow
|x|^p,~|x|^p(E_n(|x|)E^\bot_{\frac{1}{n}}(|x|)\in M_+$ and\\
$\widehat{\Phi}(|x|^pE_n(|x|)E^\bot_{\frac{1}{n}}(|x|))\leq (\sup
D(x))^p,$ we have $|x|^p\in L^1(M,\Phi)$ and\\
$\widehat{\Phi}(|x|^p)=\sup_{n\geq
1}\widehat{\Phi}(|x|^pE_n(|x|)E^\bot_{\frac{1}{n}}(|x|))$
\cite{Zak_Ch}. This means that $x\in L^p(M,\Phi)$ and
$\|x\|_{p,\Phi}\leq\sup D(x).$ Theorem \ref{3.3.} implies $\sup
D(x)\leq\|x\|_{p,\Phi},$ and therefore $\|x\|_{p,\Phi}=\sup D(x).$
\end{proof}

With the help of Theorem \ref{3.4.}, it is not difficult to show
that $L^p(M,\Phi)$ is disjointly  decomposable LNS over $S_h(
\mathcal{B})$ for all $p>1.$

\begin{thm}\label{3.5.} $(i)$ $L^p(M,\Phi)$ is a linear subspace in $S(M),$ and $\|\cdot\|_{p,\Phi}$
is the disjointly  decomposable $S_h(\mathcal{B})$-valued norm on
$L^p(M,\Phi);$

$(ii)$ $ML^p(M,\Phi)M\subset L^p(M,\Phi),$ and
$\|axb\|_{p,\Phi}\leq \|a\|_M\|b\|_M\|x\|_{p,\Phi}$ for all
$a,b\in M,~x\in L^p(M,\Phi);$

$(iii)$ If  $0\leq x\leq y\in L^p(M,\Phi),~x\in S(M),$ then $x\in
L^p(M,\Phi)$ and $\|x\|_{p,\Phi}\leq  \|y\|_{p,\Phi}.$
\end{thm}
\begin{proof} $(i)$ It is clear that $\lambda x\in L^p(M,\Phi)$ and $\|\lambda x\|_{p,\Phi}=|\lambda|\|
x\|_{p,\Phi}$ for all $x\in L^p(M,\Phi),~\lambda\in \mathbb{C}.$
Moreover, $\|x\|_{p,\Phi}\geq 0$ and $\widehat{\Phi}(|x|^p)=
\|x\|^p_{p,\Phi}=0$ if and only if $x=0.$

We claim that $x+y\in L^p(M,\Phi)$ and $\|x+y\|_{p,\Phi}\leq
\|x\|_{p,\Phi}+\|y\|_{p,\Phi}$ for each $x,y\in L^p(M,\Phi).$ By theorem \ref{3.3.}, $(x+y)z=xz+yz\in L^1(M,\Phi)$ for
all $z\in L^q(M,\Phi),$ in addition
$$
|\widehat{\Phi}((x+y)z)|\leq |\widehat{\Phi}(xz)|
+|\widehat{\Phi}(yz)|.
$$
If $\|z\|_{q,\Phi}\leq \mathbf{1}_\mathcal{B},$ then by theorem
\ref{3.4.},
$$
|\widehat{\Phi}((x+y)z)|\leq\|x\|_{p,\Phi}+\|y\|_{p,\Phi}.
$$
Using Theorem \ref{3.4.} again, we obtain  $x+y\in L^p(M,\Phi)$ and
$\|x+y\|_{p,\Phi}\leq \|x\|_{p,\Phi}+\|y\|_{p,\Phi}.$ Thus,
$L^p(M,\Phi)$ is a linear subspace in $S(M),$ and
$\|\cdot\|_{p,\Phi}$ is a $S_h(\mathcal{B})$-valued norm on
$L^p(M,\Phi).$

Let us now show that the norm $\|\cdot\|_{p,\Phi}$ is
$d$-decomposable. It is known  \cite{ChZ} that, if $x\in
L^1(M,\Phi),~\|x\|_{1,\Phi}=f_1+f_2,$ where $f_1,f_2\in
S_+(\mathcal{B}),~f_1f_2=0,$ then, setting $x_i=xp_i$ for
$p_i=\psi^{-1}(s(f_i)),~i=1,2,$ we get $x=x_1+x_2$ and
$\|x_i\|_\Phi=f_i,~i=1,2.$

Let $y\in L_+^p(M,\Phi),~\|y\|_{p,\Phi}=g_1+g_2$ where $g_1,g_2\in
S_+(\mathcal{B}),~g_1g_2=0,$ i.e. $\|y^p\|_{1,\Phi}=
\|y\|^p_{p,\Phi} =g_1^p+g_2^p.$ Set $q_i=\psi^{-1}(s(g_i^p))\in
P(\mathcal{A})\subset P(Z(M))$ and $y_i=yq_i.$ Then $ y_i^p=
y^pq_i$ and  using \cite{ChZ} for
$x=y^p,~f_i=g_i^p,~i=1,2$ we obtain that $y^pq_1+y^pq_2=y^p$ and
$\|yq_i\|_{p,\Phi}=g_i,~i=1,2.$ Since $q_1q_2=0,~q_1,q_2\in
P(Z(M)),$ we have $yq_1+yq_2=y.$

Let now $y$ be an arbitrary element from $L^p(M,\Phi)$ and let $y=u|y|$ be the polar decomposition of $y$ with the unitary
$u\in M.$ Let $\||y|\|_{p,\Phi}=\|y\|_{p,\Phi}=f_1+f_2$ where
$f_1,f_2\in S_+(\mathcal{B}),~f_1f_2=0.$ It follows from  above that for $q_i=\psi^{-1}(s(f_i^p))\in P(\mathcal{A}),$ we have
$|y|=|y|q_1+|y|q_2$ è $\||y|q_i\|_{p,\Phi}=f_i.$ Consequently,  $y=u|y|=u|y|q_1+u|y|q_2=yq_1+yq_2$ and $\|yq_i
\|_{p,\Phi}= \||yq_i|\|_{p,\Phi}=\||y|q_i\|_{p,\Phi}=f_i,~i=1,2.$
Hence, the norm $\|\cdot\|_{p,\Phi}$ is $d$-decomposable.

$(ii)$ Let  $v$ be a unitary operator in $M,$ $x\in L^p(M,
\Phi).$ Then $|vx|=(x^*v^*vx)^{\frac{1}{2}}=|x|,$ and therefore
$vx\in L^p(M,\Phi).$ Since any operator $a\in M$ is a linear
combination of four unitary operators, we have  $ax\in L^p(M,\Phi),$ due to
$(i).$

We claim that $\|ax\|_{p,\Phi}\leq \|a\|_M\|x\|_{p,\Phi}$
for $a\in M,~x\in L^p(M,\Phi).$ Let $\nu$ be a faithful normal
semifinite numerical  trace on $\mathcal{B}.$ If for some $a\in M,$
$x\in L^p(M,\Phi)$ the previous inequality is not true, then there
are $\varepsilon>0,~0\neq e\in P(\mathcal{B}),~\nu(e)<\infty$ such
that
$$e\|ax\|_{p,\Phi} \geq e\|a\|_M\|x\|_{p,\Phi}+\varepsilon e.$$
By the formula
$$
\tau(b)=\nu(e\Phi(b)(\mathbf{1}_\mathcal{B}+
\Phi(\mathbf{1})+\widehat{\Phi}(|x|^p))^{-1}), b\in Ms(\Phi_e)
$$
one defines a  faithful normal finite numerical trace on $Ms(\Phi_e).$ If $z=\psi^{-1}(e)\in P(\mathcal{A}),$ then $\Phi_e(\mathbf{1}-z)=
(\mathbf{1}_\mathcal{B}-e) e\Phi(\mathbf{1})=0,$ i.e.
$s(\Phi_e)\leq z.$ Since $\Phi(z-s(\Phi_e))=
\Phi(z(\mathbf{1}-s(\Phi_e)) =e\Phi(\mathbf{1}-s(\Phi_e))=0,$ we get
$z=s(\Phi_e).$ We  consider the $L^p$-space $L^p(Ms(\Phi_e),\tau)$
associated with the numerical trace $\tau,$ and let us show that
$xz\in L^p(Ms(\Phi_e),\tau).$ Let $x_n=E_n(|x|)|x|.$ It is clear that
$0\leq x_n^pz\uparrow |x|^pz,$ moreover
$$\tau(x_n^pz)=\nu(e\Phi(x_n^pz)(\mathbf{1}_\mathcal{B}+
\Phi(\mathbf{1})+\widehat{\Phi}(|x|^p))^{-1})\leq \nu(e)<\infty.$$
Hence, $|xz|^p=|x|^pz\in L^p(Ms(\Phi_e),\tau)$ and
$\|xz\|^p_{p,\tau}=\lim\limits_{n\to
\infty}\|x_n^pz\|^p_{p,\tau}=\nu(e\widehat{\Phi}(|x|^pz)(\mathbf{1}_\mathcal{B}+\Phi(\mathbf{1})+\widehat{\Phi}(|x|^p))^{-1}).$
Thus, if $a\in M$ then $axz\in L^p(Ms(\Phi_e),\tau),$ in
addition
$$\|a\|_M\|xz\|^p_{p,\tau}\geq
\|axz\|^p_{p,\tau}=\nu(e\widehat{\Phi}(|axz|^p)(\mathbf{1}_\mathcal{B}+\Phi(\mathbf{1})+\widehat{\Phi}(|x|^p))^{-1})=$$
$$\nu(e\|ax\|_{p,\Phi}^p)(\mathbf{1}_\mathcal{B}+\Phi(\mathbf{1})+\widehat{\Phi}(|x|^p))^{-1})\geq$$
$$ \nu(e(\|a\|_M\|x\|_{p,\Phi}+\varepsilon
)^p(\mathbf{1}_\mathcal{B}+\Phi(\mathbf{1})+\widehat{\Phi}(|x|^p))^{-1})>\|a\|_M^p\|xz\|_{p,\tau}^p,$$ which is not the case. Consequently, $\|ax\|_{p,\Phi}\leq\|a\|_M\|x\|_{p,\Phi}.$

If $b\in M,~x\in L^p(M,\Phi),$ then by  Proposition
\ref{3.2.} and from above, we have $b^*x^*\in L^p(M,\Phi).$ Using
Proposition \ref{3.2.} again, we obtain $xb=(b^*x^*)^*\in
L^p(M,\Phi)$ and $\|xb\|_{p,\Phi}=\|b^*x^*\|_{p,\Phi}\leq
\|b^*\|_M\|x^*\|_{p,\Phi}=\|b\|_M\|x\|_{p,\Phi}.$

$(iii)$ Let $0\leq x\leq y\in L^p(M,\Phi),~x\in S(M).$ It follows
from (\cite{Mur_m}, \S 2.4) that $\sqrt{x}=a\sqrt{y}$ where $a\in
M$ with $\|a\|_M\leq 1.$  Hence, $x=\sqrt{x}(\sqrt{x})^*=aya^*\in
L^p(M,\Phi)$ è $\|x\|_{p,\Phi}\leq
\|a\|_M\|a^*\|_M\|y\|_{p,\Phi}\leq \|y\|_{p,\Phi}.$

\end{proof}

Using the H\"{o}lder inequality and  the
$(bo)$-completeness of the space $(L^1(M,\Phi), \|\cdot\|_\Phi)$
we can establish the $(bo)$-completeness of  the space
$(L^p(M,\Phi), \|\cdot\|_{p,\Phi}).$

\begin{thm}\label{3.6.} Let $\Phi, M, p$ be the same as in  Theorem \ref{3.3.}. Then $(L^p(M,\Phi), \|\cdot\|_{p,\Phi})$ is
the Banach-Kantorovich space.
\end{thm}
\begin{proof} First, we assume that $\mathcal{B}$ is a $\sigma$-finite von Neumann
algebra. Then  there exists  a faithful normal finite numerical trace $\nu$ on $\mathcal{B}.$ The numerical function
$\tau(a)=\nu(\Phi(a)(\mathbf{1}_\mathcal{B}+\Phi(\mathbf{1}))^{-1})$
is a faithful normal finite trace on $M.$ Moreover,
the topology $t(M)$ coincides with
topology of convergence in measure $t_\tau$   in $(S(M),\tau)$  (\cite{Mur_m}, \S 3.5).

Let $\{x_\alpha\}_{\alpha\in A}\subset (L^p(M,\Phi), \|\cdot\|_{p,\Phi})$ be an
$(bo)$-Cauchy net i.e. $b_\gamma=\sup\limits_{\alpha,\beta\geq\gamma}\|x_\alpha
-x_\beta\|_{p,\Phi}\downarrow 0.$ According to the H\"{o}lder
inequality, for each $x\in L^p(M,\Phi)$ we have $x\in L^1(M,\Phi)$
and $\|x\|_{1,\Phi}=\widehat{\Phi}(|x|\mathbf{1})\leq
(\Phi(\mathbf{1}))^{\frac{1}{q}}\|x\|_{p,\Phi}.$ In particular,
the set $\{\|x_\alpha-x_\beta\|_{1,\Phi}
\}_{\alpha,\beta\geq\gamma}$ is bounded in $S_h(\mathcal{B}),$ and
$\sup\limits_{\alpha,\beta\geq\gamma}\|x_\alpha
-x_\beta\|_{1,\Phi}\leq(\Phi(\mathbf{1}))^{\frac{1}{q}}b_\gamma$
for all $\gamma\in A.$ Consequently  \cite{ChZ}, there exists  $x\in L^1(M,\Phi)$ such
that $\|x_\alpha -x\|_{1,\Phi}\stackrel{(o)}{\longrightarrow}0$  in particular, $x_\alpha
\stackrel{t_\tau}{\longrightarrow}x$ and  $y_\alpha=|x_\alpha
-x_\beta|\stackrel{t_\tau}{\longrightarrow}|x-x_\beta|.$ Since the function $\varphi(t)=t^p$
is continuous on $[0,\infty),$  the operator function $y
\longmapsto y^p$ is continuous on $(S_+(M),t_\tau)$ \cite{Tih}.
Hence, $0\leq y_\alpha^p \stackrel{t_\tau}{\longrightarrow}
|x-x_\beta|^p,$ in addition $\widehat{\Phi}(y_\alpha^p)=
\|x_\alpha-x_\beta\|^p_{p,\Phi}\leq b_\gamma^p.$ Using the
of Fatou's theorem \cite{Zak_Ch}, we obtain
$|x-x_\beta|^p\in L^1(M,\Phi)$ and $\widehat{\Phi}(| x-x_\beta
|^p)\leq b_\gamma^p.$ Thus, $(x-x_\beta)\in L^p(M,\Phi)$ for all
$\beta\geq\gamma$ and $\sup\limits_{\beta\geq\gamma}\|x-x_\beta
\|_{p,\Phi}\leq b_\gamma \downarrow 0.$ This means that
$x\in L^p(M,\Phi),$ and $\|x_\alpha -x\|_{p,\Phi}\stackrel{(o)}{\longrightarrow}0.$

 Now let $\mathcal{B}$ be an arbitrary  von Neumann
algebra ( not necessarily  $\sigma$-finite), and let $\{x_\alpha\}\subset L^p(M,\Phi)$ be a $(bo)$-Cauchy net.  It follows from the  above  that there exists
$x\in L^1(M,\Phi)$ such that $\|x_\alpha-x\|_{1,\Phi}
\stackrel{(o)}{\longrightarrow}0.$ In particular $x_\alpha
\stackrel{t(M)}{\longrightarrow}x.$ Let $\nu$ be a faithful normal
semifinite numerical trace on $\mathcal{B},$ and let $\{e_i\}_{i\in I}$
be the family of nonzero mutually orthogonal projections in
$\mathcal{B}$ such that $\sup\limits_{i\in I}e_i=\mathbf{1}
_\mathcal{B},$ and $\nu(e_i)<\infty$ for all $i\in I.$ It is
clear that $\{x_\alpha s(\Phi_{e_i})\}_{\alpha \in A}$ is a
$(bo)$-Cauchy net in $L^p(Ms(\Phi_{e_i}), \Phi_{e_i}).$
Since the algebra $\mathcal{B}e_i$ is $\sigma$-finite, from the above there exists  $x_i\in
L^p(Ms(\Phi_{e_i}), \Phi_{e_i})$ such that $\|x_i-x_\alpha
s(\Phi_{e_i})\|_{p,\Phi_{e_i}}\stackrel{(o)}{\longrightarrow}0.$
In particular, $x_\alpha s(\Phi_{e_i}) \stackrel{t(M
)}{\longrightarrow}x_i=x_is(\Phi_{e_i}).$ On the other hand,
convergence $x_\alpha\stackrel{t(M)}{\longrightarrow}x$ implies
$x_\alpha s(\Phi_{e_i})\stackrel{t(M)}{\longrightarrow}
xs(\Phi_{e_i}).$ Thus, $xs(\Phi_{e_i})=x_is(\Phi_{e_i})$ for all
$i\in I.$ By Proposition \ref{3.1.}, we have $x\in L^p(M,\Phi)$
and $\|x-x_\alpha\|_{p,\Phi}e_i=\|xs(\Phi_{e_i})-x_\alpha
s(\Phi_{e_i})\|_{p,\Phi_{e_i}}\stackrel{(o)}{\longrightarrow}0$
for all $i\in I$ and therefore $\|x-x_\alpha
\|_{p,\Phi}\stackrel{(o)}{\longrightarrow}0.$
\end{proof}

\begin{pr}\label{3.7.} If $\{x_\alpha\}_{\alpha\in A}\subset
L^p_h(M,\Phi)$ and $x_\alpha\downarrow 0,$ then $\|x_\alpha
\|_{p,\Phi}\downarrow 0.$
\end{pr}
\begin{proof} Let $\nu$ be a faithful normal semifinite numerical trace on $\mathcal{B}.$  If $b=\inf_{\alpha\in
I}\|x_\alpha\|_{p,\Phi}\neq 0,$ then there are $\varepsilon>0,$
$0\neq e\in P(\mathcal{B})$ with $\nu(e)<\infty$ such that
$e\|x_\alpha\|_{p,\Phi}\geq eb \geq \varepsilon e$ for all
$\alpha\in A.$ Put $\Phi_e(x)=e\Phi(x), ~x\in M,$  and $\tau(y)=\nu(\Phi(y)(\mathbf{1}_\mathcal{B}+\Phi(\mathbf{1})
+\widehat{\Phi}(x_{\alpha_0}^p))^{-1}), $  $y\in Ms(\Phi_e),$ where $\alpha_0$ is a fixed element from $A.$ Let us prove that
$L^p(Ms(\Phi_e),\tau)\subset L^p(Ms(\Phi_e),\Phi_e)$ and
$\|x\|^p_{p,\tau}=\nu(\widehat{\Phi}(|x|^p)(\mathbf{1}
_\mathcal{B}+\Phi(\mathbf{1})+\widehat{\Phi}(x^p_{\alpha_0}))^{-1})$
for all $x\in L^p(Ms(\Phi_e),\tau).$ It is sufficient to consider
the case where  $x\in L^p_+(Ms(\Phi_e),\tau).$ Set $x_n=E_n(x)
xs(\Phi_e).$ It is clear that $x_n\in (Ms(\Phi_e))_+,~x_n^p\uparrow
x^p,~x_n^p\stackrel{\tau}{\longrightarrow}x^p,$ and therefore
$x^p_n\stackrel{t(M)}{\longrightarrow}x^p.$ Moreover,
$\Phi(|x_n^p-x_m^p|)=\Phi(x^pE_n(x)E_m^\bot(x))$ as $m<n.$ Since
$\nu(e\Phi(|x_n^p-x_m^p|)(\mathbf{1}_\mathcal{B}+\Phi(\mathbf{1})
+\widehat{\Phi}(x^p_{\alpha_0}))^{-1})=\|x_n^p-x_m^p\|_{1,\tau}=
\|x^pE_n(x)E_m^\bot(x)\|_{1,\tau}\to 0$ as $n,m\to \infty,$ we get
$\Phi(|x_n^p-x_m^p|)=e\Phi(|x_n^p-x_m^p|)\stackrel{t(
\mathcal{B})}{\longrightarrow}0.$ This means that $x^p\in L^1(M,\Phi)$
and $\Phi(x^p_n)\uparrow \widehat{\Phi}(x^p),$ i.e. $x\in
L^p(Ms(\Phi_e),\Phi_e)$ è $\|x\|_{p,\Phi_e}=\sup_{n\geq
1}(\Phi(x_n^p))^{\frac{1}{p}}.$ Using the inequality $\nu(
\Phi(x_n^p)(\mathbf{1}_\mathcal{B}+\Phi(\mathbf{1})+\widehat{\Phi}(
x^p_{\alpha_0}))^{-1})=\tau(x_n^p)\leq \tau(x^p)$  we obtain that  $\widehat{\Phi}(x^p)(
\mathbf{1}_\mathcal{B}+\Phi(\mathbf{1})+\widehat{\Phi}(x^p_{\alpha_0}))^{-1}\in L_1(\mathcal{B},\nu)$
and
$$
\nu(\widehat{\Phi}(x^p)(
\mathbf{1}_\mathcal{B}+\Phi(\mathbf{1})+\widehat{\Phi}(x^p_{\alpha_0}))^{-1})=\sup_{n\geq
1}\tau(x^p_n)=\tau(x^p),
$$
i.e.
$\|x\|_{p,\tau}=\nu(\widehat{\Phi}(x^p)(\mathbf{1}_\mathcal{B}+\Phi(\mathbf{1})+\widehat{\Phi}(x^p_{\alpha_0}))^{-1}).$

Since $\{x_\alpha\}\subset L^p(M,\Phi),$ we have that $x_\alpha
s(\Phi_e)\in L^p(Ms(\Phi_e),\Phi_e),$ moreover $x_\alpha
s(\Phi_e)\downarrow 0.$ Let us show that $x=x_{\alpha_0}
s(\Phi_e)\in L^p(Ms(\Phi_e),\tau).$ As above, we consider
$x_n=E_n(x)x.$ Since
$$
0\leq \Phi(x_n^p)(\mathbf{1}_\mathcal{B}+\Phi(\mathbf{1})+
\widehat{\Phi}(x^p_{\alpha_0}))^{-1}\uparrow
\widehat{\Phi}(x^p)(\mathbf{1}_\mathcal{B}+
\Phi(\mathbf{1})+\widehat{\Phi}(x^p_{\alpha_0}))^{-1}\leq e,
$$
we get $\tau(x^p_n)\leq \nu(e)<\infty.$ Consequently,  $x\in
L^p(Ms(\Phi_e),\tau).$ The inequality $0\leq x_\alpha \leq
x_{\alpha_0},$ for $\alpha\geq\alpha_0$ implies $x_\alpha
s(\Phi_e)\in L^p(Ms(\Phi_e),\tau)$ (see Theorem
\ref{3.5.}$(iii)$). Since $x_\alpha s(\Phi_e)\downarrow 0$ and the
norm $\|\cdot\|_{p,\tau}$ is order continuous, we have $\|x_\alpha
s(\Phi_e)\|_{p,\tau}\downarrow 0,$ i.e. $\nu(e\widehat{\Phi}
(x_\alpha^p)(\mathbf{1}_\mathcal{B}+\Phi(\mathbf{1})+
\widehat{\Phi}(x^p_{\alpha_0}))^{-1})\downarrow 0.$ Hence,
$e\widehat{\Phi}(x_\alpha)^p\downarrow 0,$ which  contradicts to the
inequality $e\Phi(x_\alpha^p)\geq \varepsilon^pe.$
\end{proof}

\section{Duality for spaces $\mathbf{L^p(M,\Phi)}$ }

Let us start with the following property of $L^p$-spaces $L^p(M,\Phi).$

\begin{pr}\label{3.9.} If $x\in L^p(M,\Phi),$ $y\in L^q(M,\Phi),$
$\frac{1}{p}+\frac{1}{q}=1,$ $p,q>1,$ then $xy,yx\in L^1(M,\Phi)$
and $\widehat{\Phi}(xy)=\widehat{\Phi}(yx).$
\end{pr}
\begin{proof} Without  loss of generality, we can take $x\geq 0,~y\geq 0.$
It follows from  Theorem \ref{3.3.} that $xy\in L^1(M,\Phi).$
Hence, $yx=y^*x^*=(xy)^*\in L^1(M,\Phi)$ and
$\widehat{\Phi}(yx)=\widehat{\Phi}((xy)^*)=\overline{\widehat{\Phi}(xy)}.$
Let $x_n=xE_n(x),$ $y_n=yE_n(y).$ Then $x_n,y_n\in M_+$ and
$\|x-x_n\|_{p,\Phi}\stackrel{t(\mathcal{B})}{\longrightarrow}0,
~\|y-y_n\|_{q,\Phi}\stackrel{t(\mathcal{B})}{\longrightarrow}0.$
Using the inequalities $|\widehat{\Phi}(xy)-\Phi(x_ny_n)|\leq
|\widehat{\Phi}(xy)-\widehat{\Phi}(x_ny)|+|\widehat{\Phi}(x_ny)-\Phi(x_ny_n)|\leq
\|x-x_n\|_{p,\Phi}\|y\|_{q,\Phi}+\|x_n\|_{p,\Phi}\|y-y_n\|_{q,\Phi},$
we obtain
$\Phi(x_ny_n)\stackrel{t(\mathcal{B})}{\longrightarrow}\widehat{\Phi}(xy).$
Since $\Phi(x_ny_n)=\Phi(\sqrt{x_n}y_n\sqrt{x_n})\geq 0$ for all
$n,$ we get $\widehat{\Phi}(xy)\geq 0.$ Therefore
$\widehat{\Phi}(xy)=\overline{\widehat{\Phi}(xy)}=\widehat{\Phi}(yx).$
\end{proof}

Let $L^p(M,\Phi)^*$ be a BKS of all $S_h(\mathcal{B})$-bounded
linear mappings from $L^p(M,\Phi)$ into $S(\mathcal{B}),$ i.e.
$S_h(\mathcal{B})$ is the dual space to the BKS $L^p(M,\Phi).$ It is
clear that any $S_h(\mathcal{B})$-bounded linear operator $T$ is a
continuous mapping from $(L^p(M,\Phi),\|\cdot\|_{p,\Phi})$ into
$(S(\mathcal{B}),t(\mathcal{B})),$ i.e., if $x_\alpha,$ $x\in
L^p(M,\Phi),$ and $\|x_\alpha-x\|_{p,\Phi}
\stackrel{t(\mathcal{B})}{\longrightarrow}0,$ then $Tx_\alpha
\stackrel{t(\mathcal{B})}{\longrightarrow}Tx.$

\begin{pr}\label{3.10.} (compare with \cite{Ku1}, 5.1.9). Let $T\in L^p(M,\Phi)^*,$
$\psi:S(\mathcal{A})\to S(\mathcal{B})$ be a $*$-isomorphism from
Theorem \ref{2.4.}$(ii).$ Then $T(ax)=\psi(a)T(x)$ for all $a\in
S(\mathcal{A}),~x\in L^p(M,\Phi).$
\end{pr}
\begin{proof} By theorem \ref{2.4.}$(ii),$ for each $z\in P(\mathcal{A}),~x\in L^p(M,\Phi)$ we have $\|zx\|_{p,\Phi}=
\widehat{\Phi}(z|x|^p)^{\frac{1}{p}}=\psi(z)\widehat{\Phi}(|x|^p)^{\frac{1}{p}}=\psi(z)\|x\|_{p,\Phi}.$
Since $T\in L^p(M,\Phi)^*,$  $|Tx|\leq c\|x\|_{p,\Phi}$ for
some $c\in S_+(\mathcal{B})$ and all $x\in L^p(M,\Phi).$ Hence  $|T(zx)|\leq \psi(z)c\|x\|_{p,\Phi},$ i.e. the support
$s(T(zx))$ is majorized by the projection $\psi(z).$ Multiplying
the equality $T(x)=T(zx)+T((\mathbf{1}-z)x)$ by $\psi(z),$ we
obtain
$$\psi(z)T(x)=\psi(z)T(zx)=T(zx).$$ If $a=\sum_{i=1}^n\lambda_i
z_i$  is a simple element from $S(\mathcal{A}),$ where
$\lambda_i\in \mathbb{C},~z_i\in P(\mathcal{A}),~i=1,\dots,n,$
then
$$
T(ax)=\sum_{i=1}^n\lambda_iT(z_ix)=(\sum_{i=1}^n\lambda_i\psi(z_i))T(x)=\psi(a)T(x).
$$
Let $a$ be an arbitrary element from $S(\mathcal{A})$ and let
$\{a_n\}$ be a sequence of simple elements from $S(\mathcal{A})$
such that $a_n\stackrel{t(\mathcal{A})}{\longrightarrow}a.$ Then
$0\leq \psi(|a_n-a|)\stackrel{t(\mathcal{B})}{\longrightarrow}
0,~\psi(a_n)\stackrel{t(\mathcal{B})}{\longrightarrow}\psi(a),$
and
$$
\|a_nx-ax\|_{p,\Phi}=\widehat{\Phi}(|(a_n-a)x|^p)^{\frac{1}{p}}= \widehat{\Phi}(|a_n-a|^p|x|^p)^{\frac{1}{p}} =\psi(|a_n-a|)\|x\|_{p,\Phi}
\stackrel{t(\mathcal{B})}{\longrightarrow}0.
$$
Since $T$ is continuous,
$\psi(a_n)T(x)=T(a_nx)\stackrel{t(\mathcal{B})}{\longrightarrow}T(ax).$
Due to the convergence
$\psi(a_n)T(x)\stackrel{t(\mathcal{B})}{\longrightarrow}\psi(a)T(x),$
the proof is complete.
\end{proof}

Now we pass to description of the $S_h(\mathcal{B})$-dual space
$L^p(M,\Phi)^*.$
\begin{thm}\label{3.12.} Let $\Phi$ be an $S(\mathcal{B})$-valued Maharam trace
on the von Neumann algebra $M,$ $p,q>1,$ $\frac{1}{p}+
\frac{1}{q}=1.$

$(i)$ If $y\in L^q(M,\Phi),$ then the linear mapping
$T_y(x)=\widehat{\Phi}(xy),~x\in L^p(M,\Phi),$ is
$S(\mathcal{B})$-bounded and $\|T_y\|=\|y\|_{q,\Phi}.$

$(ii)$ If $T\in L^p(M,\Phi)^*,$ then there exists a unique $y\in
L^q(M,\Phi)$ such that $T=T_y.$
\end{thm}
\begin{proof} $(i)$ By the H\"{o}lder inequality (theorem \ref{3.3.}), $xy \in L^1(M,\Phi)$ for all
$x\in L^p(M,\Phi)$ and $|T_y(x)|= |\widehat{\Phi}(xy)|\leq
\|y\|_{q,\Phi}\|x\|_{p,\Phi}.$ Hence, $T_y$ is
$S_h(\mathcal{B})$-bounded linear mapping from $L^p(M,\Phi)$ into
$S(\mathcal{B}).$ Due to Proposition \ref{3.9.} and Theorem \ref{3.4.}
we have
$$\|T_y\|=\sup\{|\widehat{\Phi}(yx)|: x\in
L^p(M,\Phi),~\|x\|_{p,\Phi}\leq
\mathbf{1}_\mathcal{B}\}=\|y\|_{q,\Phi}.$$

$(ii)$ Since $s(\Phi(\mathbf{1}))=\mathbf{1}_\mathcal{B},$  we can
define the element $b=(\Phi(\mathbf{1}))^{-1}\in
S_+(\mathcal{A}).$ If $\Phi_1(x)=b\Phi(x),~x\in M,$ then
$L^p(M,\Phi_1)=L^p(M,\Phi)$ and
$\|x\|_{p,\Phi_1}=b^{\frac{1}{p}}\|x\|_{p,\Phi}$ for all $x\in
L^p(M,\Phi).$ Therefore, one can take
$\Phi(\mathbf{1})=\mathbf{1}_\mathcal{B}.$

Let $T\in L^p(M,\Phi)^*.$ We choose $a\in S_+(\mathcal{B})$
with $a\|T\|=s(\|T\|).$ Set $T_1(x)=aT(x),~x\in L^p(M,\Phi).$ It
is clear that $T_1\in L^p(M,\Phi)^*$ and $\|T_1\|=a\|T\|=s(\|T\|)\leq
\mathbf{1}_\mathcal{B}.$ If we show that there exists $y_1\in
L^q(M,\Phi)$ such that $T_1x=\Phi(xy_1),$ then by virtue of
Proposition \ref{3.10.}, $Tx=\|T\|T_1(xy_1)= T(x(\psi^{-1}
(\|T\|)y_1))=T(xy)$ where $y=\psi^{-1}(\|T\|)y_1\in L^q(M,\Phi).$
Thus, one can also take that $\|T\|\leq \mathbf{1}_\mathcal{B}.$

At first, we assume that the algebra $\mathcal{B}$ is  $\sigma$-finite.  Let $\nu$ be a faithful normal finite numerical trace
on $\mathcal{B}.$ Since $|\Phi(x)|\leq \|x\|_M\Phi(
\mathbf{1})\leq \|x\|_M\mathbf{1}_\mathcal{B},~x\in M,$ we get
$\Phi(x)\in L^1(\mathcal{B},\nu).$ Consider on $M$ the faithful
normal finite trace $\tau(x)=\nu(\Phi(x)),~x\in M.$ Using the same trick as
in the proof of Proposition \ref{3.7.},  we can show that
$L^p(M,\tau)\subset L^p(M,\Phi)$ and
$\tau(|x|^p)=\|x\|^p_{p,\tau}=\nu(\widehat{\Phi}(|x|^p))$ for all
$x\in L^p(M,\tau).$  Since $|T(x)|\leq \|x\|_{p,\Phi}=(\widehat{\Phi}
(|x|^p))^{\frac{1}{p}},$ we have $T(x)\in L^1(\mathcal{B},\nu)$ for
all $x\in L^p(M,\tau).$

We define on $L^p(M,\tau)$ the linear $\mathbb{C}$-valued functional
$f(x)=\nu(Tx),$ $x\in L^p(M,\tau).$ Since $|f(x)|\leq \nu(|T(x)|)\leq
\nu(\widehat{\Phi}(|x|^p)^{\frac{1}{p}}\mathbf{1}_\mathcal{B})\leq
(\nu(\widehat{\Phi}(|x|^p)))^{\frac{1}{p}}(\nu(\mathbf{1}_\mathcal{B}
))^{\frac{1}{q}}=(\nu(\mathbf{1}_\mathcal{B}))^{\frac{1}{q}}\|x\|_{p,\tau}$
for all $x\in L^p(M,\tau),$ we have that $f$ is a bounded linear
functional  on $(L^p(M,\tau),\|\cdot\|_{p,\tau}).$ Hence  there exists an operator $y\in
L^q(M,\tau)\subset L^q(M,\Phi)$ such that $f(x)=\tau(xy)$ for all
$x\in L^p(M,\tau)$ \cite{Yead2}. We claim that $\tau(xy)=\nu(\widehat{\Phi}
(xy))$ for all $x\in L^p(M,\tau).$ Let us remind that $\tau(|z|^p)=
\nu(\widehat{\Phi}(|z|^p))$ for all $z\in L^p(M,\tau).$ If $z\in
L^1_+(M,\tau),$ then $z^{\frac{1}{p}}\in L_+^p(M,\tau),$ and
therefore $\tau(z)=\nu(\widehat{\Phi}(z).$ Hence,
$\tau(z)=\nu(\widehat{\Phi}(z))$ for all $z\in L^1(M, \tau),$ in
particular, $\tau(xy)=\nu(\widehat{\Phi}(xy))$ where $x\in
L^p(M,\tau).$ Thus, $\nu(T(x))=f(x)=\tau(xy)=\nu(
\widehat{\Phi}(xy))$ for all $x\in L^p(M,\tau).$

Let $T(x)-\widehat{\Phi}(xy)=v|T(x)-\widehat{\Phi}(xy)|$ be the
polar decomposition of the element $(T(x)-\widehat{\Phi}(xy))\in
S(\mathcal{B})$ and take  $a=\psi^{-1}(v^*).$ Since
$$
0=\nu(T(ax)-\widehat{\Phi}(axy))=\nu(v^*(T(x)-\widehat{\Phi}(xy)))
=\nu(|T(x) -\widehat{\Phi}(xy)|),
$$
we have $T(x)=\widehat{\Phi}(xy)$ for all $x\in L^p(M,\tau).$

Let  $x\in L_+^p(M,\Phi), x_n=xE_n(x).$ Then
$\|x_n-x\|_{p,\Phi}\stackrel{t(\mathcal{B})}{\longrightarrow}0$ and
therefore $T(x_n)\stackrel{t(\mathcal{B})}{\longrightarrow}T(x)$
and $|\widehat{\Phi}(x_ny)-\widehat{\Phi}(xy)|\leq
\|x_n-x\|_{p,\Phi}\|y\|_{q,\Phi}\stackrel{t(\mathcal{B})}{\longrightarrow}0.$
Since $T(x_n)=\widehat{\Phi}(x_ny),$
$T(x)=\widehat{\Phi}(xy),$ i.e. $T=T_y.$

If $z$ is another element from $L^q(M,\Phi)$  with
$T(x)=\widehat{\Phi}(xz),$ $x\in L^p(M,\Phi),$ then
$\widehat{\Phi}(x(y-z))=0$ for all $x\in L^p(M,\Phi).$ Taking
$x=u^*$ where $u$ is the unitary operator from the polar
decomposition $y-z=u|y-z|,$ we obtain $\widehat{\Phi}(|y-z|)=0,$
i.e. $y=z.$

Now let $\mathcal{B}$ be a general (not necessarily  a $\sigma$-finite) von Neumann algebra.
Let $\nu$ be a faithful normal semifinite numerical trace  on
$\mathcal{B},$ and let $\{e_i\}_{i\in I}$ be a family  of
nonzero mutually orthogonal projections in $\mathcal{B}$
with $\sup\limits_{i\in I}e_i=\mathbf{1}_\mathcal{B}$ and
$\nu(e_i)<\infty$ for all $i\in I.$ It is clear that $\mathcal{B}e_i$ is a
$\sigma$-finite algebra and $\Phi_{e_i}(x)=e_i\Phi(x)$ is
$S(\mathcal{B}e_i)$-valued Maharam trace on $Ms(\Phi_{e_i}).$
Since $T\in L^p(M,\Phi)^*,$  $T_i(x)=e_iT(x)$ is
$S_h(\mathcal{B}e_i)$-bounded linear mapping onto
$L^p(Ms(\Phi_{e_i}),\Phi_{e_i}).$ By virtue of what we  proved above, there
exists the unique $y_i\in L^q(Ms(\Phi_{e_i}),\Phi_{e_i}),$ such that

$$
e_iT(xs(\Phi_{e_i}))=\widehat{\Phi_{e_i}}(xs(\Phi_{e_i})y_i)=
e_i\widehat{\Phi}(xs(\Phi_{e_i})y_i)
$$
for all $x\in L^p(M,\Phi),~i\in I.$ Moreover,
$\|y_i\|_{q,\Phi}=\|T_i\|=\|T\|e_i.$ Since $\sup\limits_{i\in
I}s(\Phi_{e_i})=\mathbf{1},$ $\{s(\Phi_{e_i})\}_{i\in I}\subset
P(Z(M)$ and $s(\Phi_{e_i})s(\Phi_{e_j})=0$ as $i\neq j,$
there exists a unique $y\in S(M)$ such that $ys(\Phi_{e_i})
=y_i.$ We have $e_i\widehat{\Phi}(|y|^q)=\widehat{\Phi}(|y_i|^q)
=\|T\|^qe_i$ for all $i\in I.$ Hence, $y\in L^q(M,\Phi)$ è
$\|y\|_{q,\Phi}=\|T\|$ (see Proposition \ref{3.1.}). In addition
$$e_i\widehat{\Phi}(xy)=\widehat{\Phi_{e_i}}(xs(\Phi_{e_i})y_i)=e_iT(xs(\Phi_{e_i}))=e_iT(x),$$
for all $i\in I,$ i.e. $T_y(x)=\widehat{\Phi}(xy)=T(x),~x\in
L^p(M,\Phi).$
\end{proof}
\begin{ls} The BKS $L^p(M,\Phi)^*$ is isometric to the
space $(L^q(M,\Phi),\|\cdot\|_{q,\Phi}).$
\end{ls}

\end{document}